\documentclass[twocolumn]{article}

\usepackage{amsmath,amsthm,amssymb}
\usepackage[hidelinks]{hyperref}
\usepackage{enumitem}
\usepackage{authblk}
\usepackage{graphicx} % For including graphics (via \includegraphics).
\usepackage{amsmath}  % Improves typographic quality of mathematical output.
\usepackage{amsfonts} % For mathematical fonts.
\usepackage{xcolor}
\usepackage{amssymb} % Extended set of mathematical symbols.
\usepackage{tikz} % Extends the figure generation capabilities of LaTeX.
\usepackage{algorithm} % Necessary for including algorithmic structures.
\newfloat{algorithm}{t}{lop}

\newcommand{\ow}{\overline{w}}
\newcommand{\oy}{\overline{Y}}
\numberwithin{equation}{section}

\newcommand{\oyy}{\overline{y}}

\newcommand{\ZZ}{\mathbb{Z}}
\newcommand{\RR}{\mathbb{R}}

\newtheorem{lemma}{Lemma}[section]
\newtheorem{theorem}[lemma]{Theorem}
\newtheorem{definition}{Definition}[section]
\newtheorem{remark}{Remark}[section]
\newtheorem{example}[remark]{Example}

\newcommand{\measurerestr}{%
  \,\raisebox{-.127ex}{\reflectbox{\rotatebox[origin=br]{-90}{$\lnot$}}}\,%
}

\newcommand{\dd}{\mathrm{d}}

\title{Geometric modelling of polycrystalline materials:\\Laguerre tessellations and periodic semi-discrete optimal transport}

\date{\today}

\author[1]{D.~P.~Bourne
}
\author[1]{M.~Pearce
}
\author[2]{S.~M.~Roper
}
\affil[1]{Maxwell Institute for Mathematical Sciences and Department of Mathematics, Heriot-Watt University, Edinburgh, UK}
\affil[2]{School of Mathematics and Statistics,
University of Glasgow, Glasgow, UK}

\begin{document}

\maketitle

\begin{abstract}
In this paper we describe a fast algorithm for generating periodic RVEs of polycrystalline materials. 
In particular, we use the damped Newton method from semi-discrete optimal transport theory to generate 3D periodic Laguerre tessellations (or power diagrams) with cells of given volumes. Complex, polydisperse RVEs with up to 100,000 grains of prescribed volumes can be created in a few minutes on a standard laptop. The damped Newton method relies on the Hessian of the objective function, which we derive by extending recent results in semi-discrete optimal transport theory to the periodic setting.   
\end{abstract}

\section{Introduction}
There is a large literature on geometric modelling of polycrystalline metals and foams using Laguerre tessellations and weighted Voronoi diagrams; see for example the following recent papers and their references: \cite{AlpersEtAl2015,BourneKokRoperSpanjer2020,KuhnEtAl2020,PerezEtAl2019,QueyRenversade2018,TeferraRowenhorst2018,Geers2021}. Applications include generating RVEs for computational homogenisation \cite{ChandrasekaranEtAl2021,GehrigEtAl2022}, fitting Laguerre tessellations to imaging data of polycrystalline microstructures \cite{PetrichEtAl2021,PetrichEtAl2019}, and modelling grain growth \cite{AlpersEtAl2022}. 

This paper builds on the research programme initiated 
%by the authors and Tata Steel Research \& Development 
in \cite{BourneKokRoperSpanjer2020}, where recent results from optimal transport theory \cite{Santambrogio2015} were exploited to develop fast algorithms for generating Laguerre tessellations with grains of given volumes. These ideas were %picked up in the engineering community and 
developed further by \cite{KuhnEtAl2020} (see below) and applied by \cite{ChandrasekaranEtAl2021} to study biopolymer aerogels. In this paper we improve the speed of the algorithm from \cite{BourneKokRoperSpanjer2020}.

This paper also extends some theoretical results in semi-discrete optimal transport theory from \cite{KitagawaMerigotThibert2019,MerigotThibert2020} to the periodic quadratic cost function (see Theorem \ref{thm:RegK}). This extension to periodic domains was driven by applications, not only by the application to microstructure modelling (where the RVEs should be periodic to avoid artificial boundary effects in computational homogenisation), but also by a recent application in weather modelling \cite{EganEtAl2022}.

\paragraph{Outline of the paper} Section \ref{section:maths} includes the mathematical theory of periodic semi-discrete optimal transport. In particular, in Section \ref{subsec:OT problem} we recall that periodic
Laguerre tessellations with cells of given volumes can 
generated by maximising the concave function $\mathcal{K}^\Lambda$, defined in equation \eqref{eq:Kant}.
%, which is equivalent to solving the equation $\nabla \mathcal{K}^\Lambda=0$.
In \cite{BourneKokRoperSpanjer2020} $\mathcal{K}^\Lambda$ is maximised using a 1st-order method (the 2nd-order damped Newton method of \cite{KitagawaMerigotThibert2019} was also proposed but not implemented). Faster optimisation methods for maximising $\mathcal{K}^\Lambda$ were implemented in \cite{KuhnEtAl2020}, including the Barzilai-Borwein method and a modified Newton method. 

In this paper we implement the damped Newton method of \cite{KitagawaMerigotThibert2019}, where the concave function $\mathcal{K}^\Lambda$ is maximised by applying Newton's method with a bespoke backtracking scheme to the nonlinear equation $\nabla \mathcal{K}^\Lambda=0$. This requires an expression for the Hessian of  $\mathcal{K}^\Lambda$, which does not appear in the literature as far as we are aware (\cite{KitagawaMerigotThibert2019} does not include it since the periodic quadratic transport cost is not differentiable). We derive it in Theorem \ref{thm:RegK}. 
%Section \ref{Subsec:Reg} can largely be skipped by practitioners only interested in the application to microstructure modelling, with the exception of Theorem \ref{thm:RegK}.

In Section \ref{section:damped Newton} we state the damped Newton method of \cite{KitagawaMerigotThibert2019} before applying it to microstructure modelling in Section \ref{section:RVEs}. In particular, we combine the damped Newton method with \cite[Algorithm 2]{BourneKokRoperSpanjer2020} to develop a fast algorithm for generating RVEs with grains of given volumes. In Example \ref{example:log normal} we generate a polydisperse RVE with 10,000 grains of given volumes in less than a minute on a standard laptop, and in Section \ref{section:Numerical tests} we generate RVEs with 100,000 grains in a matter of minutes.
In Section \ref{section:Numerical tests} we also study the 
%run time of the damped Newton method and the 
number of Newton iterations and backtracking steps for the damped Newton method. 
%We show that we can generate RVEs with 100,000 grains in a matter of minutes.

\paragraph{Summary of main contributions} 
\begin{itemize}[leftmargin=*]
\item Theorem \ref{thm:RegK}: We prove that the Kantorovich function $\mathcal{K}^\Lambda$ is twice differentiable (under suitable assumptions) and compute its Hessian. This extension of some results from \cite{deGouray2019,KitagawaMerigotThibert2019,MerigotThibert2020} to the periodic setting, while relatively straightforward, is important for applications.

\item Knowledge exchange: We show how the damped Newton method from the mathematical theory of semi-discrete optimal transport \cite{KitagawaMerigotThibert2019} can be used to generate RVEs of polycrystalline materials with grains of given volumes.
\item Software: All our code is available on GitHub.
\smallskip
\\
\emph{MATLAB-Voro}: MATLAB mex files for generating 2D and 3D periodic and non-periodic Laguerre tessellations using Voro++ \cite{Voro++}.
\\ \url{https://github.com/smr29git/MATLAB-Voro}
\medskip
\\
\emph{MATLAB-SDOT}: MATLAB functions for solving semi-discrete optimal transport problems  using the damped Newton method.
\\
\url{https://github.com/DPBourne/MATLAB-SDOT}
\medskip
\\
\emph{Laguerre-Polycrystalline-Microstructures}: MATLAB functions for generating RVEs of polycrystalline microstructures using Laguerre tessellations, including all the examples from this paper.
\url{https://github.com/DPBourne/Laguerre-Polycrystalline-Microstructures}
\end{itemize}

%\paragraph{Outline of the paper} Theoretical results in Section 2 (can largely be skipped by practitioners except for statement of Hessian). Damped Newton method of \cite{KitagawaMerigotThibert2019} recalled in Section 3. Application to geometric modelling in Section 4, where we generate a periodic RVE .... Numerical study of damped Newton in Section 5 (run times, show how few Newton and back-tracking steps are needed).

%Motivation: Driven by applications. The need to extend the theory of damped Newton method to the periodic setting came from applications in the steel industry (periodic RVEs) and weather modelling (work with Charlie).

%Novelty: Damped Newton method doesn't apply directly to the periodic setting. We show how it can be easily extending to this case. 

%Make code available on GitHub as a selling point of the paper.

\section{Periodic semi-discrete optimal transport}
\label{section:maths}

\subsection{Notation}
Throughout this paper $|x|$ denotes the standard Euclidean norm of a vector $x \in \mathbb{R}^d$.

\paragraph{Periodic domain} First we define a general class of periodic domains that includes the flat cylinder, flat torus and triply-periodic cuboid.
Let $\Lambda \subset \RR^d$, $d \ge 2$, be the integer span of $1\leq k \leq d$ linearly independent vectors $v_1,...,v_k \in \RR^d$, i.e., 
\[
\Lambda = \text{span}_{\ZZ}\{v_1,...,v_k\}.
\]
If $k=d$, then $\Lambda$ is a \emph{lattice}. Let 
\begin{align*}
%\label{defn:VcellLambda}
V = \{x \in \RR^d:|x| \le |x-u| \; \forall \; u \in \Lambda  \}.
    \end{align*}
In other words, $V$ is the fundamental \emph{Voronoi cell} in the periodic Voronoi tessellation of $\mathbb{R}^d$ generated by $\Lambda$. For example, if $d=2$, $k=1$, $v_1=(1,0)$, then $V=[-1/2,1/2]\times \mathbb{R}$ is the fundamental domain of the flat cyclinder. If     $d=2$, $k=2$, $v_1=(1,0)$, $v_2=(0,1)$, then $V=[-1/2,1/2]\times [-1/2,1/2]$ is the fundamental domain of the flat torus. For the microstructure application in Section \ref{section:RVEs}, $\Lambda$ and $V$ are given in equations \eqref{eq:Lambda} and \eqref{eq:V}.
    
\paragraph{Transport cost} Let $c_\Lambda:\mathbb{R}^d \times \mathbb{R}^d \to \mathbb{R}$ be the periodic quadratic transport cost defined by
 \begin{align*}
    %\label{defn:PeriodicCost}
     c_{\Lambda}(x,y)=|x-y|^{2}_{\Lambda}:=\min_{u \in \Lambda}|x-y-u|^2.
 \end{align*}
 Observe that $c_\Lambda$ is continuous but not differentiable, which is why the results of \cite{deGouray2019,KitagawaMerigotThibert2019,MerigotThibert2020} do not immediately apply. For example, for the case of the flat torus given above with $d=2$, $k=2$, $v_1=(1,0)$, $v_2=(0,1)$, then $c_{\Lambda}((0,0),(t,0))=\min \{ t^2,(1-t)^2\}$ for $t \in [0,1]$, which is not differentiable as a function of $t$ at $t=1/2$.

\paragraph{Source measure} 
Let $\rho \in L^1(V;[0,\infty))$ and let $\mu = \rho \, \mathcal{L}^d \measurerestr V$ be the measure on $V$ that is 
absolutely continuous with respect to the Lebesgue measure with density $\rho$. For example, for the microstructure application we take $d=3$, $\rho(x)=1$ for all $x \in V$ so that $\mu(A)$ is the volume of a set $A \subseteq V$. In general, we assume that the support of $\mu$ is contained in a compact, convex subset $X$ of $V$ and that the restriction of $\rho$ to $X$ is continuous. 
 
 Extend $\rho$ to $\mathbb{R}^d$ by $\Lambda$-periodicity, namely, define $\rho_\Lambda:\mathbb{R}^d \to [0,\infty)$ by $\rho_\Lambda(x+u)=\rho(x)$ for all $x \in V$, $u \in \Lambda$ (this defines $\rho_\Lambda$ uniquely a.e.). Define $\mu_\Lambda = \rho_\Lambda \,  \mathcal{L}^d$. Then $\mu_\Lambda$ satisfies
%$\mu_{\Lambda} \in \mathcal{M}(\RR^d)$ satisfies 
$\mu_\Lambda \measurerestr V = \mu$ and 
    \begin{align}
        \label{eq:MuExtension}
        \mu_{\Lambda}(A) = \mu_{\Lambda}(A+u) 
    \end{align} 
for all $u \in \Lambda$ and Lebesgue-measurable sets $A \subset \mathbb{R}^d$. 
%Let $\mu \in \mathcal{M}(V)$ have compact support 
%and be absolutely continuous with respect to the Lebesgue measure. Let $\mu_{\Lambda} \in \mathcal{M}(\RR^d)$ be the $\Lambda$-periodic extension of $\mu$ satisfying $\mu_\Lambda \measurerestr V = \mu$ and 
%    \begin{align}
%        \label{eq:MuExtension}
%        \mu_{\Lambda}(A) = \mu_{\Lambda}(A+u) 
%    \end{align} 
%for all $u \in \Lambda$ and all Lebesgue-measurable sets $A \subset \mathbb{R}^d$. 
%To be precise, let $\mu = \rho \mathcal{L}^d$ with $\rho \in C_c(V;[0,\infty))$. {\color{red}[Too strong. We just want $\rho$ to be continuous on its support, whic may be a subset of $V$.]} Extend $\rho$ to $\mathbb{R}^d$ by $\Lambda$-periodicity, namely, define $\rho_\Lambda:\mathbb{R}^d \to [0,\infty)$ by $\rho_\Lambda(x+u)=\rho(x)$ for all $x \in V$, $u \in \Lambda$ (this defines $\rho_\Lambda$ uniquely almost everywhere). Then define $\mu_\Lambda = \rho_\Lambda \mathcal{L}^d$. For the regularity results in Section \ref{Subsec:Reg} we need to assume that $\rho_\Lambda$ is continuous, which essentially means that $\rho$ is continuous and satisfies `periodic boundary conditions'.  {\color{red}[Check that we need this assumption. Compare our notation with M\'erigot's: What is $\Omega_X$, $X$, etc. $\Omega_X=\mathbb{R}^d$, $X \subseteq V$, need $\rho$ continuous on $X$.]}

\paragraph{Target measure} Let 
%$Y=\{ y_1,\ldots,y_n\} \subset V$ 
$Y \subset \mathrm{int}(V)$
be a finite set. We call its elements \emph{seeds}. 
%In an abuse of notion, we will sometimes write $Y=(y_1,\ldots,y_n) \in V^n$. 
Let $m:Y \to (0,\infty)$, $y \mapsto m_y$. Define the discrete measure $\nu := \sum_{y \in Y} m_y \delta_y$.
%\in \mathcal{M}(V)$.
We assume that $\sum_{y \in Y} m_y = \mu(V)$ so that $\nu(V)=\mu(V)$. Throughout this paper $n=\# Y=\# \mathrm{spt}(\nu)$. For the microstructure application the $m_y$ are the target volumes of the $n$ grains.

\paragraph{Periodic Laguerre tessellations}
For each seed $y \in Y$, we associate a \emph{weight} $w_y \in \mathbb{R}$. Let $w$ be the \emph{weight map} $w:Y \to \mathbb{R}$, $y \mapsto w_y$.
For each $y \in Y$, we define the \emph{periodic Laguerre cell}
 \begin{align*}
 & L^{\Lambda}_{y}(w;Y) := %L^{c_\Lambda}_{y}(w;V,Y) 
 \\
 & \{x \in V: |x-y|^{2}_{\Lambda}-w_y \leq |x-z|^{2}_{\Lambda}-w_z\; \forall \; z \in Y\}.
 \end{align*}
 The collection of all cells $\{  L^{\Lambda}_{y}(w;Y) \}_{y \in Y}$ is 
 the \emph{periodic Laguerre tessellation of $V$} generated by $(Y,w)$. 
 %\sr{Perhaps \emph{periodic Laguerre tessellation of $V$} since the notation includes $\Lambda$ and then the language refers to that aspect of the notation too}
 
%\paragraph{Periodic semi-discrete optimal transport}
%Given a Borel map $T:V \to Y$, define its Monge cost by
%\[
%M^\Lambda(T) =\int_V |T(x)-x|_\Lambda^2 \, \dd \mu(x).
%\]
%\[
%M^\Lambda(T) = \sum_{y \in Y} \int_{T^{-1}(\{y\})} %|y-x|_\Lambda^2 \, \dd \mu(x).
%\]
\subsection{The optimal transport problem}
\label{subsec:OT problem}
The \emph{Kantorovich function} $\mathcal{K}^\Lambda:\mathbb{R}^{n} \to \mathbb{R}$ is
\begin{align}
\nonumber
        \mathcal{K}^\Lambda(w) & = \sum_{y \in Y} \int_{L^{\Lambda}_{y}(w;Y)} (|x-y|^2_\Lambda-w_y) \, \text{d}\mu
        (x) \\
        \label{eq:Kant}
        & 
        \quad 
        + \sum_{y \in Y}m_y w_y.
    \end{align}   
The (dual)  periodic semi-discrete optimal transport problem is the optimisation problem
\[
\max \left\{ \mathcal{K}^\Lambda(w) : w \in \mathbb{R}^n \right\}.
\]
The map
$\mathcal{K}^\Lambda$ is concave and its critical points satisfy
\[
\mu \big(L^\Lambda_{y}(w;Y)\big) = m_y
\]
(see Theorem \ref{thm:RegK}).
Hence, if $w = \mathrm{argmax} \; \mathcal{K}^\Lambda$,
then the periodic Laguerre tessellation  $\{  L^{\Lambda}_{y}(w;Y) \}_{y \in Y}$ has cells of  masses $\{m_y\}_{y \in Y}$. For example, the cells have volumes $\{m_y\}_{y \in Y}$ if $\rho=1$.
%Moreover, for completeness we recall that the  map $T:V \to Y$,
%\[
%T(x)=\mathrm{argmin}_{y \in Y} (c_\Lambda(x,y)-w_y),
%\]
%is the optimal transport map for the (primal) optimisation problem of transporting $\mu$ to $\nu$ with periodic quadratic
%cost $c_\Lambda$ \cite{Santambrogio2015}. 
%\db{We might have to cut the remark about $T$ to save space. We don't use $T$ anywhere in the paper. This remark could go in Mason's thesis.}
%{\color{red}[Cut this remark to save space? I personally think it rounds things off nicely! - Mason]}

\subsection{Regularity of the Kantorovich function}
\label{Subsec:Reg}
In this section we prove that $\mathcal{K}^\Lambda$ is twice continuously differentiable and compute its first and second derivatives. To overcome the lack of smoothness of $c^\Lambda$, we rewrite $L^{\Lambda}_{y}$ and $\mathcal{K}^\Lambda$ in terms of the standard quadratic cost and standard Laguerre cells. Then the regularity of $\mathcal{K}^\Lambda$ follows easily from  \cite{KitagawaMerigotThibert2019,MerigotThibert2020}.

   \begin{definition}%[Extended set of generators]
        \label{defn:EquivRelation}
        Given $y_1,y_2 \in \RR^d$, we write $y_1 \sim y_2$ if $y_1-y_2 \in \Lambda$. 
        We say that $\overline{y} \in \mathbb{R}^d$ is a 
        periodic copy of a seed $y \in Y$ if $\overline{y} \sim y$.
        Define $\overline{Y}$ to be the set of all periodic copies of the seeds in $Y$, namely
        \[
        \overline{Y} = Y + \Lambda = \{ y + u : y \in Y, \; u \in \Lambda \}.
        \]
        If $\oyy \in \oy$,
         $\overline{y} \sim y \in Y$, we assign $\oyy$ the weight $w_{\oyy}:=w_y$. 
         Note that if $y_1 \sim y_2$, then $w_{y_1} = w_{y_2}$.
         Define the extended weight map $\ow: \overline{Y} \to \mathbb{R}$, $\oyy \mapsto w_{\oyy}$.
    Given a domain $\Omega \subseteq \mathbb{R}^d$, the non-periodic Laguerre tessellation $\{L_{\oyy}(\ow;\Omega,\oy)\}_{\oyy \in \oy}$ of $\Omega$ generated by the extended set of seeds and weights $(\oy,\ow)$ is given by
\begin{align*}
& L_{\oyy}(\ow;\Omega,\oy) := 
\\
& \{x \in \Omega: |x-\oyy|^{2}-w_{\oyy} \leq |x-z|^{2}-w_z\; \forall \; z \in \oy \}.
\end{align*}
%for all $\oyy \in \oy$.
\end{definition}
\begin{lemma}
\label{Lemma:3.1}
    The periodic Laguerre cells can be written in terms of the non-periodic cells as 
        \begin{equation}
        \label{eq:Lemma3.1}
        L_y^{\Lambda}(w;Y) = \bigcup_{u \in \Lambda} L_{y+u}(\ow;V,\oy).
        \end{equation}
\end{lemma}

\begin{proof}
Let $x \in L_y^{\Lambda}(w;Y)$. By definition we have
\[ 
|x-y|^2_{\Lambda}-w_y \leq |x-z|^2_{\Lambda}-w_z \quad \forall \; z \in Y.
\]
Hence there exists $u_{xy} \in \Lambda$ such that,
for all $z \in Y$,
\begin{align*}
|x-y-u_{xy}|^2 -w_y 
& \le |x-z|^2_{\Lambda}-w_z 
\\
& \leq |x-z-u|^2-w_z \quad \forall \; u \in \Lambda.
\end{align*}
Observe that $w_y=w_{y+u_{xy}}$ and $w_z=w_{z+u}$.
Therefore the previous inequality can be restated as 
\[ 
|x-(y+u_{xy})|^2 -w_{y+u_{xy}} \leq |x-\overline{z}|^2-w_{\overline{z}} \quad \forall \; \overline{z} \in \oy.
\]
Therefore $x \in L_{y+u_{xy}}(\ow;V,\oy)$ and so
\[
L_y^{\Lambda}(w;Y) \subseteq \bigcup_{u \in \Lambda} L_{y+u}(\ow;V,\oy).
\]

Now we prove the reverse inclusion. 
Let $x \in L_{y+u}(\ow;V,\oy)$ for some $u \in \Lambda$. By definition,
%of $\| \cdot \|_\Lambda$, %and $\overline{w}$,
\begin{align*}
|x-y|_{\Lambda}^2 - w_y 
& \le |x-y-u|^2 - w_{y+u}
\\
& \le |x-\overline{z}|^2-w_{\overline{z}} \quad \forall \; \overline{z} \in \oy
\end{align*}
where the second inequality follows from the fact that $x \in L_{y+u}(\ow;V,\oy)$. Therefore
\[
|x-y|_{\Lambda}^2 - w_y 
\le |x-z-u|^2-w_z \quad \forall \; z \in Y, \, u \in \Lambda.
\]
Taking the minimum over $u$  gives
\[
|x-y|_{\Lambda}^2 - w_y 
\le |x-z|_{\Lambda}^2-w_z \quad \forall \; z \in Y.
\]
Therefore $x \in L_y^{\Lambda}(w;Y) $, as required. 
%$\square$
\end{proof}

%{\color{red}[Add a figure to illustrate Lemmas 2.1 \& 2.2? Maybe not enough space in the special issue.]}

\begin{lemma}
\label{Lemma:3.2}
    The $\mu$-measure of the periodic Laguerre cells can be expressed in terms of the measure of the non-periodic Laguerre cells as follows:
        \begin{align*}
        \mu(L^{\Lambda}_{y}(w;Y))
        & =
        \sum_{u \in \Lambda} \mu \left( 
        L_{y+u}(\overline{w};V,\overline{Y})
        \right) 
        \\
        & =
        \mu_\Lambda(L_y(\overline{w};\RR^d,\overline{Y})).
        \end{align*}
\end{lemma}

\begin{proof}
The first equality follows immediately from Lemma \ref{Lemma:3.1} and the fact that the Laguerre cells 
$\{ L_{y+u}(\overline{w};V,\overline{Y}) \}_{u \in \Lambda}$ are disjoint up to a set of measure zero.
Now we turn our attention to the second inequality. We claim that
\begin{equation}
\label{eq:Morrissey}
L_y(\overline{w};\RR^d,\overline{Y}) 
= \bigcup_{u \in \Lambda} \left( L_{y+u}(\overline{w};V,\overline{Y}) - u \right).
\end{equation}
Let $x \in L_y(\overline{w};\RR^d,\overline{Y})$. Then 
\[ 
|x-y|^2-w_y \leq |x-z|^2-w_z \quad \forall \; z \in \overline{Y}.
\]
Choose $u \in \Lambda$ so that $x + u \in V$. Then we can rewrite the previous inequality as
\[
|x+u-(y+u)|^2-w_{y+u} \leq |x-z|^2-w_z \quad \forall\; z \in \overline{Y}.
\]
Therefore $x+u \in L_{y+u}(\overline{w};V,\overline{Y})$ and so 
\[ 
L_y(\overline{w};\RR^d,\overline{Y}) \subseteq \bigcup_{u \in \Lambda }
\left(L_{y+u}(\overline{w};V,\overline{Y}) - u \right).
\]
Conversely, take $x \in L_{y+u}(\overline{w};V,\overline{Y}) - u $ for some $u \in \Lambda$. For all $z \in \overline{Y}$,
\[ 
|x-y|^2 -w_y
= |x+u-(y+u)|^2-w_y 
\leq |x-z|^2-w_z.%\;\; \forall \;z \in \overline{Y}. 
\]
Therefore $x \in L_y(\overline{w};\RR^d,\overline{Y}) $ and \eqref{eq:Morrissey} follows. 

Next we prove that $\{ L_{y+u}(\overline{w};V,\overline{Y}) - u \}_{u \in \Lambda}$ are disjoint sets up to a set of measure zero. 
Suppose that $x \in L_{y+u_i}(\overline{w};V,\overline{Y})-u_i$ for $i \in \{1,2\}$,
%and $x \in L_{y+u_2}(\overline{w};V,\overline{Y})-u_2$ for
$u_1,u_2 \in \Lambda$, $u_1 \neq u_2$. Then
$x+u_i \in V$. % and $x+u_2 \in V$. 
By definition of $V$,
\[
|x+u_i| \le |x+u_i-u| \quad \forall \; u \in \Lambda
\]
for $i \in \{1,2\}$. That is,
\[
|x-(-u_i)| \le |x-u| \quad \forall \; u \in \Lambda.
\]
In other words, $x$ lies in the Voronoi cells with generators $-u_1$ and $-u_2$ in the Voronoi tessellation of $\mathbb{R}^d$ generated by $\Lambda$. But the intersection of Voronoi cells is a set of measure zero, as desired. 

Combining everything and using $\mu_{\Lambda} \measurerestr V=\mu$ gives %and $\mu_\Lambda$ is $\Lambda$-periodic gives
\begin{align*}
& \sum_{u \in \Lambda} \mu \left( 
        L_{y+u}(\overline{w};V,\overline{Y})
        \right)
\\
&         
\stackrel{\eqref{eq:MuExtension}}{=} 
        \sum_{u \in \Lambda} \mu_{\Lambda} \left( 
        L_{y+u}(\overline{w};V,\overline{Y}) - u
        \right)
\\
& 
\stackrel{\eqref{eq:Morrissey}}{=} \mu_{\Lambda} \left( L_y(\overline{w};\RR^d,\overline{Y})  \right)
\end{align*}
as required. %$\square$
\end{proof}

\begin{lemma}
\label{Lemma:3.3}
$\mathcal{K}^\Lambda$
%The Kantorovich function %$\mathcal{K}^\Lambda$ 
can be written in terms of the non-periodic transport cost and Laguerre cells as
\begin{align}
\nonumber
& \mathcal{K}^\Lambda (w) - \sum_{y \in Y} m_y w_y
\\
\label{eq:Kv1}
& =\sum_{y \in Y}\int_{L_y(\ow;\RR^{d},\oy)}\left(|x-y|^2 - w_y \right) \, \dd \mu_\Lambda (x)
\\
\label{eq:Kv2}
& = \sum_{y \in \oy}\int_{L_y(\ow;V,\oy)}\left(|x-y|^2 - w_y \right) \, \dd \mu (x).
\end{align}
\end{lemma}

%We will use \eqref{eq:Kv1} for implementation and \eqref{eq:Kv2} for establishing the regularity of $\mathcal{K}^\Lambda$.

\begin{proof}
Let $y \in Y$, $u \in \Lambda$, $x \in L_{y+u}(\ow;V,\oy)$. We claim that
\begin{equation}
\label{eq:ProofLemma3.3_1}
|x-y|_{\Lambda} = 
|x-y-u|.
\end{equation}
To prove this observe that
%If $x \in L_{y+u}(\ow;V,\oy)$, then
\begin{align*}
|x-y-u|^2 
& = |x-(y+u)|^2 - w_{y+u} + w_{y+u}
\\
& \le 
|x-(y+v)|^2 \underbrace{- w_{y+v} + w_{y+u}}_{=0}% \qquad \forall \; v \in \Lambda.
\end{align*}
for all $v \in \Lambda$. 
Therefore
\[
|x-y-u| = \min_{v \in \Lambda} |x-y-v| = |x-y|_{\Lambda}
\]
as claimed. 
For all $y \in Y$,
%\begin{align*}
%& \int_{L^{\Lambda}_{y}(w;Y)} |x-y|^2_\Lambda \, \text{d}\mu(x)
%\\
%& = 
%\int_{\bigcup_{u \in \Lambda} L_{y+u}(\ow;V,\oy)}
%|x-y|^2_\Lambda \, \text{d}\mu(x) 
%& (\textrm{by Lemma \ref{Lemma:3.1}})
%\\
%& =
%\int_{\bigcup_{u \in \Lambda} L_{y+u}(\ow;V,\oy)}
%|x-y-u|^2 \, \text{d}\mu(x) 
%& (\textrm{by \eqref{eq:ProofLemma3.2_1}})
%\\
%& = 
%\int_{\bigcup_{u \in \Lambda} (L_{y+u}(\ow;V,\oy)-u)}
%|x'-y|^2 \, \text{d}\mu_\Lambda(x') 
%& (x'=x-u, \textrm{ \eqref{eq:MuExtension}})
%\\
%& = 
%\int_{L_y(\overline{w};\RR^d,\overline{Y}) }
%|x'-y|^2 \, \text{d}\mu_\Lambda(x') 
%& (\textrm{by \eqref{eq:Morrissey}})
%\end{align*}
%Combining this result with Lemma \ref{Lemma:3.2} and the definition of $\mathcal{K}^\Lambda$ completes the proof.
\begin{align}
\nonumber
& \int_{L^{\Lambda}_{y}(w;Y)} (|x-y|^2_\Lambda - w_y) \, \text{d}\mu(x)
\\
\nonumber
& \stackrel{\eqref{eq:Lemma3.1}}{=} 
\int_{\bigcup_{u \in \Lambda} L_{y+u}(\ow;V,\oy)}
(|x-y|^2_\Lambda - w_y) \, \text{d}\mu(x) 
%& (\textrm{by Lemma \ref{Lemma:3.1}})
\\
\label{eq:ProofLemma3.3_2}
& \stackrel{\eqref{eq:ProofLemma3.3_1}}{=}
\int_{\bigcup_{u \in \Lambda} L_{y+u}(\ow;V,\oy)}
(|x-y-u|^2 - w_y) \, \text{d}\mu(x) 
%& (\textrm{by \eqref{eq:ProofLemma3.3_1}})
%\end{align}
%Therefore
%\begin{align}
%\nonumber
%& \sum_{y \in Y} \int_{L^{\Lambda}_{y}(w;Y)} (|x-y|^2_\Lambda %- w_y) \, \text{d}\mu(x)
\\
\nonumber
& =\sum_{u \in \Lambda} \int_{L_{y+u}(\ow;V,\oy)}
(|x-(y+u)|^2 - w_y) \, \text{d}\mu(x).
%& (\textrm{by \eqref{eq:ProofLemma3.3_2}})
\end{align}
%\label{eq:ProofLemma3.3_3}
Combining this with the definition of $\mathcal{K}^\Lambda(w)$ gives %we have
\begin{align}
\nonumber
&\mathcal{K}^\Lambda(w)-\sum_{y\in Y}m_yw_y=\\
\nonumber
&=\sum_{y\in Y}\sum_{u \in \Lambda} \int_{L_{y+u}(\ow;V,\oy)}
(|x-(y+u)|^2 - w_y) \, \text{d}\mu(x)
\\
\nonumber
& = 
\sum_{\overline{y} \in \oy} \int_{L_{\overline{y}}(\ow;V,\oy)}
(|x-\overline{y}|^2 - w_{\overline{y}}) \, \text{d}\mu(x). \end{align}
%since $w_{y}=w_{y+u}$ for all $y \in Y$, $u \in \Lambda$.
%Combining %\eqref{eq:ProofLemma3.3_3} 
%with 
%the definition of $\mathcal{K}^\Lambda$ proves \eqref{eq:Kv2}.
This proves \eqref{eq:Kv2}.

Using the substitution $x'=x-u$ in \eqref{eq:ProofLemma3.3_2}  gives
%Alternatively, we can rewrite equation \eqref{eq:ProofLemma3.3_2} as
\begin{align*}
& \int_{L^{\Lambda}_{y}(w;Y)} (|x-y|^2_\Lambda - w_y) \, \text{d}\mu(x)
\\
& = 
\int_{\bigcup_{u \in \Lambda} (L_{y+u}(\ow;V,\oy)-u)}
(|x'-y|^2 - w_y) \, \text{d}\mu_\Lambda(x') 
%& (x'=x-u, \textrm{ \eqref{eq:MuExtension}})
\\
& = 
\int_{L_y(\overline{w};\RR^d,\overline{Y}) }
(|x'-y|^2 - w_y) \, \text{d}\mu_\Lambda(x') 
%& (\textrm{by \eqref{eq:Morrissey}})
\end{align*}
by \eqref{eq:Morrissey}.
This 
%Combining this with the definition of $\mathcal{K}^\Lambda$ 
proves \eqref{eq:Kv1}, as required. %$\square$
%and completes the proof.
\end{proof}

In Theorems \ref{thm:RegQ} and \ref{thm:RegK} we enumerate 
the seeds $y_1,\ldots,y_M$, $M \in \mathbb{N}$, and in an abuse of notation let $w$ denote both the weight map $w:y_i \mapsto w_{y_i} \in \mathbb{R}$ and the vector $(w_1,\ldots,w_M):=(w_{y_1},\ldots,w_{y_M})$.

%In the following two theorems we take 
%$Y=\{y_1,\ldots,y_n\} \subset V$ and, in an abuse of notation,  let $w$ denote both the map $w:Y \to \mathbb{R}$, $y_i \mapsto w_{y_i}$,  and the vector $(w_1,\ldots,w_n):=(w_{y_1},\ldots,w_{y_n})\in\mathbb{R}^n$. 

First we recall a result from
\cite{KitagawaMerigotThibert2019} (Theorems 1.3 \& 4.1) about the regularity of the Kantorovich function for the standard (non-periodic) quadratic cost (see also \cite[Proposition 2]{deGouray2019}, \cite[Theorem 45]{MerigotThibert2020}). 

\begin{theorem}[\cite{KitagawaMerigotThibert2019}]
\label{thm:RegQ}
Let $Y=\{y_1,\ldots,y_N\} \subset V$ be a set of distinct seeds.
Define 
$g:\mathbb{R}^{N} \to \mathbb{R}$  by
\begin{align*}
        g(w) = \sum_{y \in Y} \int_{L_{y}(w;V,Y)} (|x-y|^2-w_y) \, \dd \mu(x). % + \sum_{\tilde{y} \in \tilde{Y}}m_{\tilde{y}} w_{\tilde{y}}.
    \end{align*} 
%\begin{align*}
%        g(w) = \sum_{i=1}^n \int_{L_{y_i}(w;V,Y)} (|x-y_i|^2-w_i) \, \dd \mu(x). % + \sum_{\tilde{y} \in \tilde{Y}}m_{\tilde{y}} w_{\tilde{y}}.
%    \end{align*} 
Then $g \in C^1(\mathbb{R}^N)$ is concave with 
\[ 
\frac{\partial g}{\partial w_i}(w) = %m_{\tilde{y}_i} 
- \mu \big(L_{y_i}(w;V,Y)\big), \quad i \in \{1,\ldots,N\}. 
\]
For $i,j \in \{1,\ldots,N\}$, define
\[
L_{y_i y_j}(w;V,Y) = L_{y_i}(w;V,Y) \cap L_{y_j}(w;V,Y)
\]
(this may be the empty set).
On the set 
\begin{equation}
\label{eq:MassAssump}
\left\{ w \in \mathbb{R}^N : \mu \big( L_{y}(w;V,Y) \big) > 0 \; \forall \, y \in Y \right\} 
\end{equation}
$g$ is twice continuously differentiable and, for $i \ne j$, %with Hessian
 \begin{align*}
 \frac{\partial^2 g}{\partial w_{i} \partial w_{j}}(w)
 & =  \int_{L_{y_i y_j}(w;V,Y)} \frac{\rho(x)}{2|y_i-y_j|} \, \dd \mathcal{H}^{d-1}(x), 
% \; i \ne j,
 %\quad \forall \; 
% i,j \in \{ 1,\ldots,N\}, \; 
% i \neq j,
\\     
\frac{\partial^2 g}{\partial w_i^2}(w) & = - \sum_{\substack{j=1 \\ j \ne i}}^N \frac{\partial^2 g}{\partial w_{i} \partial w_{j}}(w).
\end{align*}
\end{theorem}

%\begin{proof}
%{\color{red}See \cite[Theorems 1.1, 1.3, 4.1]{KitagawaMerigotThibert2019}, where this result is proved for a general class of cost functions on Riemannian manifolds. See also \cite[Theorems 40, 45]{MerigotThibert2020} for a simplified proof in Euclidean space. We use the opposite sign convention for the weights $w_i$ to \cite{KitagawaMerigotThibert2019,MerigotThibert2020}. The assumption in \cite[Theorem 45 \& Lemma 46]{MerigotThibert2020} that the set of seeds is \emph{generic} with respect to the cost and $\partial X$ is not needed here since we assume that $X$ is convex and $\mu ( L_{y}(w;V,Y) ) > 0$ (which ensure that the set $S$ in the proof of \cite[Lemma 46]{MerigotThibert2020} has zero $(d-1)$-Hausdorff measure).} $\square$
%\end{proof}

Now we extend Theorem \ref{thm:RegQ} to the periodic cost.

%By combing Theorem \ref{thm:RegQ} with the lemmas above we obtain the following.
\begin{theorem}%[Regularity of $\mathcal{K}^\Lambda$]
\label{thm:RegK}
Let $Y=\{y_1,\ldots,y_n\} \subset \mathrm{int}(V)$ be a set of distinct seeds.
The 
Kantorovich function $\mathcal{K}^\Lambda \in C^1(\mathbb{R}^n)$ is concave with
%$\mathcal{K}^\Lambda \in C^1(\mathbb{R}^n)$ is concave  
%with
\[ 
\frac{\partial \mathcal{K}^\Lambda}{\partial w_i}(w) = m_{y_i} - \mu \big(L^\Lambda_{y_i}(w;Y)\big), \quad i \in \{1,\ldots,n\}. 
\]
Fix $w \in \mathbb{R}^n$. Let %$\tilde{Y}$ be the set
\[
\tilde{Y} = \tilde{Y}(w) = \{ y \in \oy : L_y(\ow;V,\oy) \ne \emptyset \} 
\]
and let $\tilde{w}=\ow|_{\tilde{Y}}$ be the restriction of
$\ow$ to $\tilde{Y}$.
If 
\begin{equation}
\label{eq:PositiveVolumeAssumption}
\mu \big( L_y(\tilde{w};V,\tilde{Y}) \big) > 0 \quad \forall \; y \in \tilde{Y},
\end{equation}
then $\mathcal{K}^\Lambda$ is twice continuously differentiable at $w$ and, for $i \ne j$, 
 \begin{align*}
 \frac{\partial^2 \mathcal{K}^\Lambda}{\partial w_{i} \partial w_{j}}(w)
 & =  
 \sum_{\substack{y,y' \in \oy\\y \sim y_i\\y' \sim y_j}} 
\int_{L_{yy'}(\ow;V,\oy)} \frac{\rho(x)}{2|y-y'|} \, \dd \mathcal{H}^{d-1}
\\     
\frac{\partial^2 \mathcal{K}^\Lambda}{\partial w_i^2}(w) & = - \sum_{\substack{j=1 \\ j \ne i}}^n \frac{\partial^2 \mathcal{K}^\Lambda}{\partial w_{i} \partial w_{j}}(w).
\end{align*}
\end{theorem}

\begin{proof}
First note that the set $\tilde{Y}$ is finite since $V$ has finitely many neighbours. %{\color{red}[Add proof in thesis.]}

Step 1: First we prove that for all $y \in \tilde{Y}$,
\begin{equation}
\label{eq:Proof Theorem 2.5 0}    
L_y(\tilde{w};V,\tilde{Y})=L_y(\ow;V,\oy).
\end{equation}
If $x \in L_y(\ow;V,\oy)$, then
\[
| x-y |^2 - w_y \le |x-z|^2-w_z
\]
for all $z \in \oy$, hence for all $z \in \tilde{Y}$ since $\tilde{Y} \subseteq \oy$. Therefore $L_y(\ow;V,\oy) \subseteq$
$L_y(\tilde{w};V,\tilde{Y})$. Now we prove the opposite inclusion. Let $x \in L_y(\tilde{w};V,\tilde{Y})$. Suppose for contradiction that $x \notin L_y(\ow;V,\oy)$. Then $x \in L_z(\ow;V,\oy)$ for some $z \in \oy$, $z \ne y$. Moreover, $z \in \tilde{Y}$ else $L_z(\ow;V,\oy)=\emptyset$. Therefore $x \in L_y(\tilde{w};V,\tilde{Y}) \cap L_z(\tilde{w};V,\tilde{Y})$ and so
\begin{equation}
\label{eq:Proof Theorem 2.5 1}
|x-y|^2 - w_y = |x-z|^2 - w_z.
\end{equation}
Since $x \notin L_y(\ow;V,\oy)$, there exists $u \in \oy$ such that
\begin{equation}
\label{eq:Proof Theorem 2.5 2}
|x-y|^2 - w_y > |x-u|^2 - w_u.
\end{equation}
Since $x \in L_z(\ow;V,\oy)$,
\begin{equation}
\label{eq:Proof Theorem 2.5 3}
|x-z|^2-w_z \le |x-v|^2 - w_v \quad \forall \; v \in \oy.
\end{equation}
Combining \eqref{eq:Proof Theorem 2.5 1}-\eqref{eq:Proof Theorem 2.5 3} gives
\begin{align*}
|x-z|^2-w_z & = |x-y|^2 - w_y 
 \\ & > |x-u|^2 - w_u 
 \ge |x-z|^2-w_z,
\end{align*}
which contradicts the assumption  $x \notin L_y(\ow;V,\oy)$.

Step 2: Fix $w \in \mathbb{R}^n$. We prove that, under assumption \eqref{eq:PositiveVolumeAssumption},  $\tilde{Y}(w)$ is independent of $w$ in a neighbourhood of $w$. Let
\[
    \varepsilon = \frac 12 
\inf_{y \in \oy \setminus \tilde{Y}} \inf_{x \in V} 
f(x,y)
\]
where $f:V\times \oy \rightarrow\mathbb{R}$ is given by
\[
f(x,y) = |x-y|^2 - w_y - \min_{z \in \tilde{Y}} (|x-z|^2-w_z).
\]
Let $x \in V$, $y \in \oy \setminus \tilde{Y}$. Then $x \in L_u(\tilde{w};V,\tilde{Y})$ for some $u \in \tilde{Y}$, $u \ne y$ and so
\[
f(x,y)=|x-y|^2 - w_y - (|x-u|^2 - w_u)>0.
\]
Both infimums in the definition of $\varepsilon$ are  minimums because $f(\cdot,y)$ is a continuous function on the compact set $V$ (since $f(\cdot,y)$ is the pointwise minimum of a family of uniformly Lipschitz functions) and because we can restrict the infimum over $\oy \setminus \tilde{Y}$ in the definition of $\varepsilon$ to a finite set (since $f(x,y)$ blows up as $|y|\to \infty$). Therefore $\varepsilon > 0$.

Let $\psi:Y \to \mathbb{R}$, $|(\psi_{y_1},\ldots,\psi_{y_n})| < \varepsilon$, $x \in V$. Assume that $y \in \oy$, $y \notin \tilde{Y}(w)$. 
Then $x \in L_u(\tilde{w};V,\tilde{Y})$ for some $u \ne y$, $u \in \tilde{Y}$. 
Therefore
%We show that $x \notin L_y(\tilde{w}+\tilde{\psi};V,\tilde{Y})$ where $\tilde{\psi}=\overline{\psi}|_{\tilde{Y}}$. We have
\begin{align*}
&|x-y|^2-(w_y+\psi_y) - [|x-u|^2-(w_u+\psi_u)]    
\\
& = |x-y|^2 - w_y - \min_{z \in \tilde{Y}} [|x-z|^2-w_z] + \psi_u - \psi_y
\\
& > 
|x-y|^2 - w_y - \min_{z \in \tilde{Y}} [|x-z|^2-w_z] - 2 \varepsilon \; \ge 0.
%\\
%& \ge 0.
\end{align*}
Therefore $x \notin L_y(\overline{w}+\overline{\psi};V,\oy)$.  Since $x \in V$ was arbitrary, this implies that
$L_y(\overline{w}+\overline{\psi};V,\oy)=\emptyset$. Hence 
$y \notin \tilde{Y}(w+\psi)$, and $\oy  \setminus \tilde{Y}(w) \subseteq \oy \setminus \tilde{Y}(w+\psi)$ for all $|\psi|<\varepsilon$. 

To complete the proof of Step 2 we need to show that, under assumption \eqref{eq:PositiveVolumeAssumption}, $\tilde{Y}(w) \subseteq \tilde{Y}(w+\psi)$ if $\psi$ is sufficiently small. Let $y \in \tilde{Y}(w)$.
The map $w \mapsto \mu(L_y(\overline{w};V,\oy))$ is continuous (cf.~\cite[Proposition 38(vii)]{MerigotThibert2020}). Therefore there exists $\delta_y>0$ such that 
$\mu(L_y(\ow+\overline{\psi};V,\oy))>0$ if $|\psi|<\delta_y$, which implies that
$y \in \tilde{Y}(w+\psi)$.
Take $\delta = \min_{y \in \tilde{Y}(w)} \delta_y > 0$ since $\tilde{Y}(w)$ is finite. Therefore 
$\tilde{Y}(w) \subseteq \tilde{Y}(w+\psi)$ for all 
$|\psi|<\delta$, as required.

Step 3: Now we compute the gradient and Hessian of $\mathcal{K}^\Lambda$. Fix $w_0 \in \mathbb{R}^n$. By 
\eqref{eq:Kv2} and \eqref{eq:Proof Theorem 2.5 0}, for $w$ in a sufficiently small neighbourhood of $w_0$,
\[
\mathcal{K}^\Lambda(w)=g(\tilde{w}) + \sum_{y \in Y} m_y w_y,
\]
where $g$ is as in Theorem \ref{thm:RegQ} with
$N=|\tilde{Y}(w_0)| \ge n$ and with the set of seeds $\tilde{Y}(w_0)$ (which is fixed in the chosen neighbourhood of $w_0$). 

Observe that $\mathcal{K}^\Lambda$ is concave since it is the sum of a linear function and the composition of the concave function $g$ with the linear function $w \mapsto \tilde{w}$. 

The gradient of $\mathcal{K}^\Lambda$ follows from
the Chain Rule and Theorem \ref{thm:RegQ}:
\begin{align*}
\frac{\partial \mathcal{K}^\Lambda}{\partial w_y}(w_0)
& = \sum_{z \in \tilde{Y}(w_0)} \frac{\partial g}{\partial \tilde{w}_z}(\tilde{w}_0) \frac{\partial \tilde{w}_z}{\partial w_y}(w_0)
+ m_y
\\
& = 
\sum_{\substack{z \in \tilde{Y}(w_0)\\z \sim y}} \frac{\partial g}{\partial \tilde{w}_z}(\tilde{w}_0)
+ m_y
\\
& = 
-\sum_{\substack{z \in \tilde{Y}(w_0)\\z \sim y}} 
\mu \big(L_{z}(\tilde{w}_0;V,\tilde{Y}(w_0))\big)
+ m_y
\\
& =- \mu(L^{\Lambda}_{y}(w_0;Y)) + m_y 
\end{align*}
by Lemma \ref{Lemma:3.2} and equation \eqref{eq:Proof Theorem 2.5 0}. 

If assumption \eqref{eq:PositiveVolumeAssumption} holds for $w=w_0$, then $\tilde{w}_0$ belongs to the set \eqref{eq:MassAssump} (with seeds $\tilde{Y}(w_0))$. Therefore $w \mapsto g(\tilde{w})$ is twice continuously differentiable at $w_0$ and so is $\mathcal{K}^\Lambda$. Let $y_i,y_j \in Y$, $i \ne j$. Applying the Chain Rule and Theorem \ref{thm:RegQ} gives 
\begin{align*}
& \frac{\partial^2 \mathcal{K}^\Lambda}{\partial w_i \partial w_j} (w_0)
\\
& = \sum_{\substack{z \in \tilde{Y}(w_0)\\z \sim y_i}} 
\sum_{z' \in \tilde{Y}(w_0)}
\frac{\partial^2 g}{\partial \tilde{w}_z \partial \tilde{w}_{z'}}(\tilde{w}_0)
\frac{\partial \tilde{w}_{z'}}{\partial w_{y_j}}(w_0)
\\
& = \sum_{\substack{z \in \tilde{Y}(w_0)\\z \sim y_i}} 
\sum_{\substack{z' \in \tilde{Y}(w_0)\\z' \sim y_j}} 
\frac{\partial^2 g}{\partial \tilde{w}_z \partial \tilde{w}_{z'}}(\tilde{w}_0)
\\
& = \sum_{\substack{z,z' \in \tilde{Y}(w_0)\\z \sim y_i, z' \sim y_j}} 
\int_{L_{zz'}(\tilde{w}_0;V,\tilde{Y}(w_0))} \frac{\rho(x)}{2|z-z'|} \, \dd \mathcal{H}^{d-1}(x)
\\
& = \sum_{\substack{z,z' \in \oy\\z \sim y_i, z' \sim y_j}} 
\int_{L_{zz'}(\ow_0;V,\oy)} \frac{\rho(x)}{2|z-z'|} \, \dd \mathcal{H}^{d-1}(x)
\end{align*}
by equation \eqref{eq:Proof Theorem 2.5 0}. We also used the fact that $z \ne z'$ if $z \sim y_i$, $z'\sim y_j$, $y_i \ne y_j$, which follows from the assumption $Y \subset \mathrm{int}(V)$.

Finally, the diagonal entries of the Hessian are obtained by differentiating the following expression with respect to $w_i$:
\[
\sum_{j=1}^n \frac{\partial \mathcal{K}^\Lambda}{\partial w_j}(w_0)
= \sum_{j=1}^n 
\left[ m_{y_j} - \mu \big(L^\Lambda_{y_j}(w_0;Y)\big) \right]
= 0.
\]
This concludes the proof. %$\square$
\end{proof}

\begin{remark}
The Hessian of $\mathcal{K}^\Lambda$ can be rewritten as follows, which we found more convenient for computational purposes: For $i \ne j$,
\[
 \frac{\partial^2 \mathcal{K}^\Lambda}{\partial w_{i} \partial w_{j}}(w) = \sum_{\substack{y \in \oy\\y \sim y_j}}
 \int_{L_{y_i y}(\ow;\mathbb{R}^d,\oy)} \frac{\rho(x)}{2|y_i-y|} \, \dd \mathcal{H}^{d-1}.
\]
%{\color{red}[Include proof in thesis.]}
\end{remark}

\section{The damped Newton method} 
\label{section:damped Newton}
Given $w \in \mathbb{R}^n$, let $H(w) \in \mathbb{R}^{n \times n}$ denote the Hessian matrix $D^2 \mathcal{K}(w)$.
If $w$ satisfies  \eqref{eq:PositiveVolumeAssumption}, it can be shown that $H(w)$ is singular with 1-dimensional kernel spanned by $(1,1,\ldots,1) \in \mathbb{R}^n$ ($H(w)$ is a weighted graph Laplacian matrix of a connected graph). Moreover, the $(n-1)\times(n-1)$ matrix $\hat{H}(w)$ obtained by deleting the last row and column of $H(w)$ is non-singular. Let $e(w)$ be the error
%the error in the mass of the Laguerre cells:
\[
e(w) = 
\max_{y \in Y} \left| \frac{\partial \mathcal{K}^\Lambda}{\partial w_y}(w) \right|
=
\max_{y \in Y} 
\left\| \mu \big(L^\Lambda_{y}(w;Y)\big) -m_y\right\|. 
\]

We recall the damped Newton method of Kitagawa, M\'erigot and Thibert \cite{KitagawaMerigotThibert2019} for maximising $\mathcal{K}^\Lambda$:
\paragraph{Initialisation} Choose $w^{0} \in \mathbb{R}^n$ such that
%all the Laguerre cells $L^\Lambda_{y}(w^{0};Y)$ have positive mass, i.e., 
\begin{equation}
\label{eq:epsilon}
    \mu \big(L^\Lambda_{y}(w^{0};Y)\big) > 0 \quad \forall \; y \in Y,
\end{equation}
i.e., all the Laguerre cells have positive mass. Let
\[
\varepsilon := \frac 12 \min 
\left\{ 
\min_{y \in Y} \mu \big(L^\Lambda_{y}(w^{0};Y)\big) , \, \min_{y \in Y} m_y
\right\} >0.
\]

\paragraph{Iteration step} Give $w^{k-1}$, define $w^{k}$ as follows:
\begin{enumerate}
    \item Define the Newton direction $d^{k} \in \mathbb{R}^{n-1}$ by solving the sparse SPD linear system
    \[
    -\hat{H}(w^{k-1}) \, d^{k} = b
    \]
    where $b_i= \partial \mathcal{K}^\Lambda / \partial w_i (w^{k-1})$, $i =1,\ldots,n-1$. 
    %Define $d^{(k+1)} \in \mathbb{R}^n$ by $d^{(k+1)}_i=\hat{d}^{(k+1)}_i$,
    %$i =1,\ldots,n-1$, $d^{(k+1)}_n=0$.
    \item
    Find the smallest value of $l_k \in \{0\} \cup \mathbb{N}$ such that $w^{k,l_k} \in \mathbb{R}^n$ defined by
    $w^{k,l_k}_n =0$, 
    \begin{align*}
    w^{k,l_k}_i  =w^{k-1}_i+2^{-l_k}d^{k}_i, \, i \in \{1,\ldots,n-1\},
    \end{align*}
    satisfies 
    \begin{gather*}
    \min_{y \in Y} \mu \big(L^\Lambda_{y}(w^{k,l_k};Y)\big) \ge \epsilon,
    \\
    e(w^{k,l_k}) \le (1-2^{-(l_k+1)}) \, e(w^{k-1}).
    \end{gather*}
    \item Define the Newton update $w^{k}:=w^{k,l_k}$.
    \end{enumerate}

\paragraph{Stopping condition} Terminate the algorithm when the mass percentage error is less than some prescribed tolerance $\eta>0$:
\[
100 \cdot \max_{y \in Y} \frac{\left| \mu \big(L^\Lambda_{y}(w^{k};Y)\big) -m_y\right|}{m_y} < \eta.
\]

\smallskip

We refer to $l_k$ as the number of \emph{backtracking steps} at iteration $k$. If $l_k=0$, then $w^{k}$ is the standard Newton step for the nonlinear system $\nabla \mathcal{K}^\Lambda(w)=0$. The backtracking ensures that the iterates $w^{k}$ remain in the region where $\mathcal{K}^\Lambda$ is twice differentiable.

In \cite{KitagawaMerigotThibert2019} it is proved that, for a class of continuously differentiable transport costs $c$, the damped Newton method converges for any admissible initial guess with linear rate and asymptotic quadratic rate (see \cite[Prop.~6.1]{KitagawaMerigotThibert2019}). While their result does not apply to our non-smooth cost $c_\Lambda$, we have observed quadratic convergence in  numerical experiments and believe the proof can be easily extended.

\section{RVEs of polycrystalline materials} 
\label{section:RVEs}
We can apply the theory above to generate 3D polycrystalline microstructures.
%by 3D triply-periodic Laguerre tessellations with grains of given volumes. 
We take $d=3$, 
\begin{gather}
\label{eq:Lambda}
\Lambda=\mathrm{span}_\mathbb{Z}\{(L_1,0,0),(0,L_2,0),(0,0,L_3)\},
\\
\label{eq:V}
V=\left[-\tfrac{L_1}{2},\tfrac{L_1}{2} \right] \times \left[-\tfrac{L_2}{2},\tfrac{L_2}{2} \right] \times
\left[-\tfrac{L_3}{2},\tfrac{L_3}{2} \right], 
\end{gather}
where $L_1,L_2,L_3>0$,
and $\rho(x)=1$ for all $x \in V$. 
Then $\mu(L_y^\Lambda(w;V))$ equals the volume of $L_y^\Lambda(w;V)$.
We can generate a regularised periodic Laguerre tessellation of $V$ with grains of given volumes by combining
\cite[Algorithm 2]{BourneKokRoperSpanjer2020} with the damped Newton method as follows:

\paragraph{Input}
The number of grains $n$, the desired grain volumes $m_1,\ldots,m_n$ with $\sum_i m_i=L_1L_2L_3$, the volume percentage error tolerance $\eta>0$, and the number of regularisation (Lloyd) steps $K$.

\paragraph{Initialisation}
Randomly select $n$ distinct seeds $y^{(0)}_1,\ldots,y^{(0)}_n \in \mathrm{int}(V)$. 
Set $Y_{(0)}:=\{ y^{(0)}_1,\ldots,y^{(0)}_n \}$,
$w_{(0)}:=0 \in \mathbb{R}^n$, and
\[
L_i^{(0)}:=L_{y^{(0)}_i}(\overline{w_{(0)}};\mathbb{R}^d,\overline{Y_{(0)}}), \quad i \in \{1,\ldots,n\}.
\]

\paragraph{Iteration step} For $k=1,\ldots,K$: \begin{enumerate}
    \item \emph{Regularisation}. Define $y^{(k)}_i$ to be the centroid 
    \[
    y_i^{(k)} = \frac{1}{\mu(L_i^{(k-1)})} \int_{L_i^{(k-1)}} x \, \dd x.
    \]
    Set $Y_{(k)}=\{ y^{(k)}_1,\ldots,y^{(k)}_n \}$. 
    \item \emph{Damped Newton}. Apply the damped Newton method with initial guess $w^0=0 \in \mathbb{R}^n$ to find a weight vector $w_{(k)} \in \mathbb{R}^n$ satisfying
    \[
100 \cdot \max_{i} \frac{\Big| \mu \Big(L^\Lambda_{y_i^{(k)}}(w_{(k)};Y_{(k)})\Big) -m_i\Big|}{m_i} < \eta.
\]
For $i \in \{1,\ldots,n\}$, define
\[
L_i^{(k)}:=L_{y^{(k)}_i}(\overline{w_{(k)}};\mathbb{R}^d,\overline{Y_{(k)}}).
\]
\end{enumerate}

\paragraph{Output} Define $y_i:=y_i^{(K)}$, $Y=\{y_1,\ldots,y_n\}$, and $w=w_{(K)}$. Then the periodic Laguerre cells
$L^\Lambda_{y_1}(w;V),\ldots,L^\Lambda_{y_n}(w;V)$ have volumes $m_1,\ldots,m_n$ up to $\eta$ percentage error.

\medskip

If $K$ is large, then the outputted Laguerre tessellation is approximately a \emph{centroidal Laguerre tessellation}
\cite{BourneRoper2015}. This means that the seeds $y_i$ are approximately the centroids of the grains, which results in grains that are more `regular' (round). 

The initial guess $w^0=0$ for the damped Newton steps satisfies assumption \eqref{eq:epsilon} since the corresponding Laguerre tessellation is a Voronoi tessellation, and Voronoi cells have positive volume when the seeds lie in the box.

The main difference between this algorithm and the one implemented in \cite{BourneKokRoperSpanjer2020} is that we compute the weights using the 2nd-order damped Newton method instead of the slower 1st-order BFGS method. A rigorous study of various optimisation algorithms (steepest descent, Malitsky-Mishchenko, Barzilai-Borwein, BFGS, modified Newton) and regularisation steps (Lloyd, Anderson-accelerated Lloyd) is given in \cite{KuhnEtAl2020}.

\begin{example}
\label{example:log normal}
We reproduce Example 5.5 from \cite{BourneKokRoperSpanjer2020}.
We take $n=10,000$ grains, $L_1=L_2=L_3=2$, $K=5$, $\eta=1$.
The target grain volumes are drawn from a log-normal distribution exactly as in \cite[Example 5.5]{BourneKokRoperSpanjer2020}. The RVE is shown in Figure \ref{fig:Example_4_1}. This example took 47.24\,s on a laptop with processor Intel(R) Core(TM) i5-1135G7 @ 2.40GHz (cf.~the run time of 669 s in \cite[Example 5.5]{BourneKokRoperSpanjer2020}). 
\begin{figure}[h!]
    \centering
    \includegraphics[trim={0 0.1cm 0 0},clip,width=\linewidth]{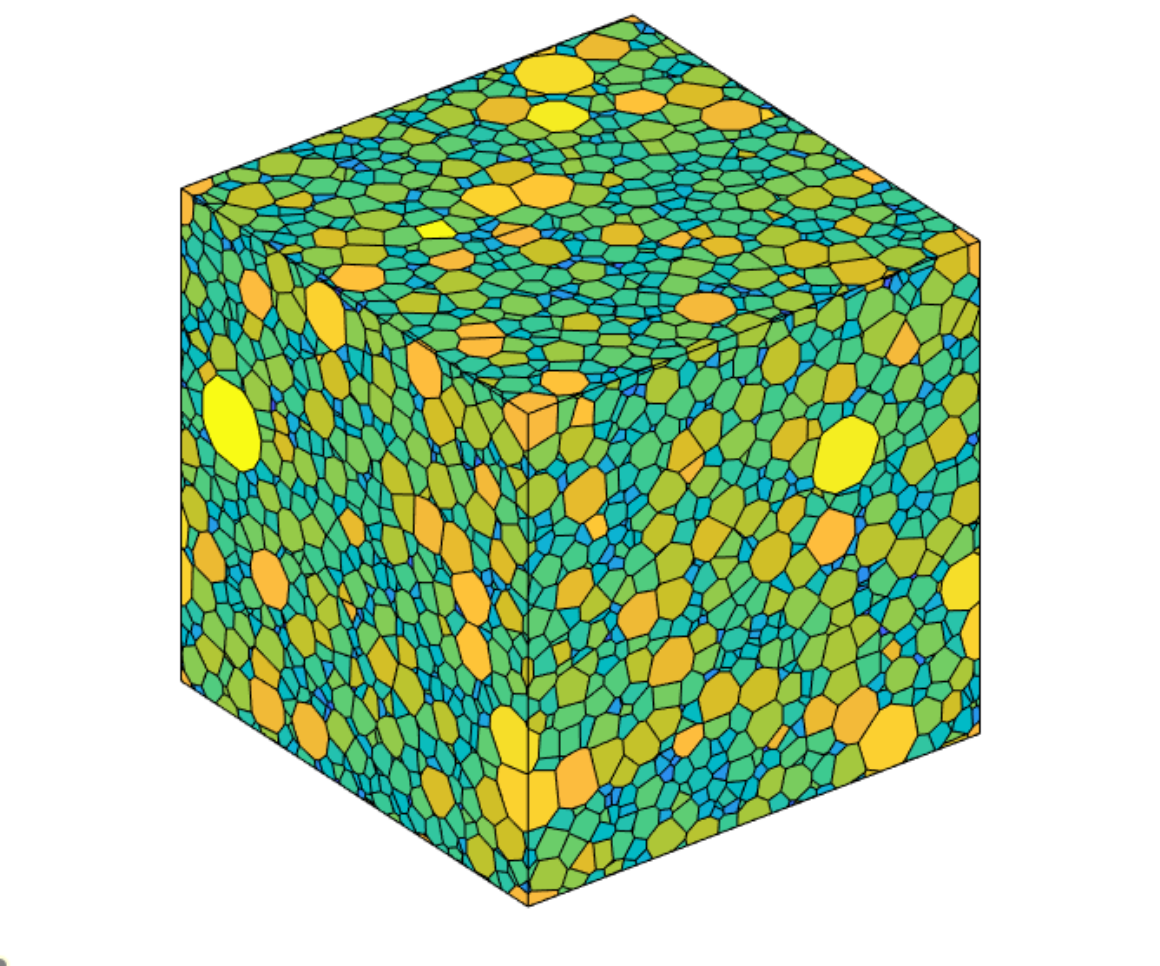}
    \caption{A periodic Laguerre tessellation with $10,000$ grains of prescribed volumes (up to $1$\% error), where the target volumes were drawn from a log-normal distribution; see Example \ref{example:log normal}. The grains are coloured according to their volume using a log scale.}
    \label{fig:Example_4_1}
\end{figure}
\end{example}
 
\section{Numerical tests of damped Newton}
\label{section:Numerical tests}
In this section we illustrate the performance of the damped Newton method from Section \ref{section:damped Newton}. In both examples $\Lambda$ and $V$ are given by \eqref{eq:Lambda}, \eqref{eq:V}
with $L_1=L_2=L_3=1$ and $\rho(x)=1$ for all $x \in V$ so that the damped Newton method generates Laguerre cells of given volumes.

\begin{example}
\label{example:run times}
Table \ref{table:run times} reports mean run times of the damped Newton method over 100 numerical experiments for three types of microstructure.

The 1st column is the number of grains. The 2nd column is the mean run time for an idealised single phase (SP) microstructure, where all the grains have the same volume $m_i=L_1L_2L_3/n$ for all $i$.
The 3rd column is the mean run time for an idealised dual phase (DP) microstructure, where half of the grains have volume $x$ and the other half have volume $5x$, where $x$ is such that $\sum_i m_i = L_1 L_2 L_3$. The 4th column (log-normal) is the mean run time for a more realistic microstructure where the grain volumes are drawn from a log-normal distribution as in Example \ref{example:log normal}.
In all three experiments the seeds $y_1,\ldots,y_n$ were drawn randomly from a uniform distribution and we average the run times over 100 draws. The volume tolerance is $\eta=1\%$.
We used the same laptop that was used for Example \ref{example:log normal}. 

The run times grow superlinearly with respect to $n$ but subquadratically. They increase as the complexity of the microstructure increases. 
%A direct comparison with the run times in  \cite[Table 2 (SP), Table 5 (DP), Table 6 (log-normal)]{KuhnEtAl2020} is not easy since the largest value of $n$ is $16,000$ in
%\cite[Tables 2 \& 5]{KuhnEtAl2020} and $1000$ in 
%\cite[Table 6]{KuhnEtAl2020}, the processors are different, and the seed locations are different (a single random draw of the seeds is used in \cite{KuhnEtAl2020}). Our implementation of the damped Newton method appears to be faster than the modified Newton method used in \cite{KuhnEtAl2020} but slower than their Barzilai–Borwein method, at least for single phase and dual phase RVEs with  $n = 16,000$. {\color{red}[Remove this comparison to \cite{KuhnEtAl2020} since it isn't very scientific?]}

\begin{table}[h!]
\centering
\begin{tabular}{c|c|c|c}
\hline
$n$ & \multicolumn{3}{c}{mean run time (s)}
\\
& SP & DP & log-normal
\\
\hline
100 &  0.03 & 0.09 & 0.06 \\
250 & 0.07 & 0.16 & 0.23 \\
500 & 0.18 & 0.34 & 0.53 \\
1000 & 0.48 & 0.69 & 1.29 \\
2500 & 1.30 & 1.82 & 3.51 \\
5000 & 2.83 & 4.29 & 7.97 \\
10,000 & 6.68 & 9.11 & 18.15 \\
25,000 & 22.22 & 29.26 & 57.84 \\
50,000 & 57.60 & 72.53 & 224.63 \\
100,000 & 171.65 & 222.42 & 396.26 \\
\hline
\end{tabular}
\caption{Run times of the damped Newton method for generating Laguerre tessellations
with $n$ grains of given volumes (see Example \ref{example:run times}), where all the grains have the same volume (SP), there are two grain sizes (DP), or the volumes are drawn from a log-normal distribution (log-normal).}
\label{table:run times}
\end{table}
\end{example}
     
\begin{example}     
\label{example:back tracking}
Figure \ref{fig:backtrack} shows 
the number of Newton iterations and the number of backtracking steps $l_k$ per Newton iteration for 100 random log-normal microstructures with $n=100,000$ grains (as in the previous example). 
%The volume tolerance is $\eta=1\%$.
There is no backtracking at all after the 4th Newton iteration in any of the 100 experiments. In most experiments there is only one backtracking step in the 3rd Newton iteration.
We do not observe the extensive backtracking in the early iterations reported in \cite[Section 3.4]{KuhnEtAl2020}. 
\begin{figure}[h!]
\begin{center}
\includegraphics[clip,width=\linewidth]{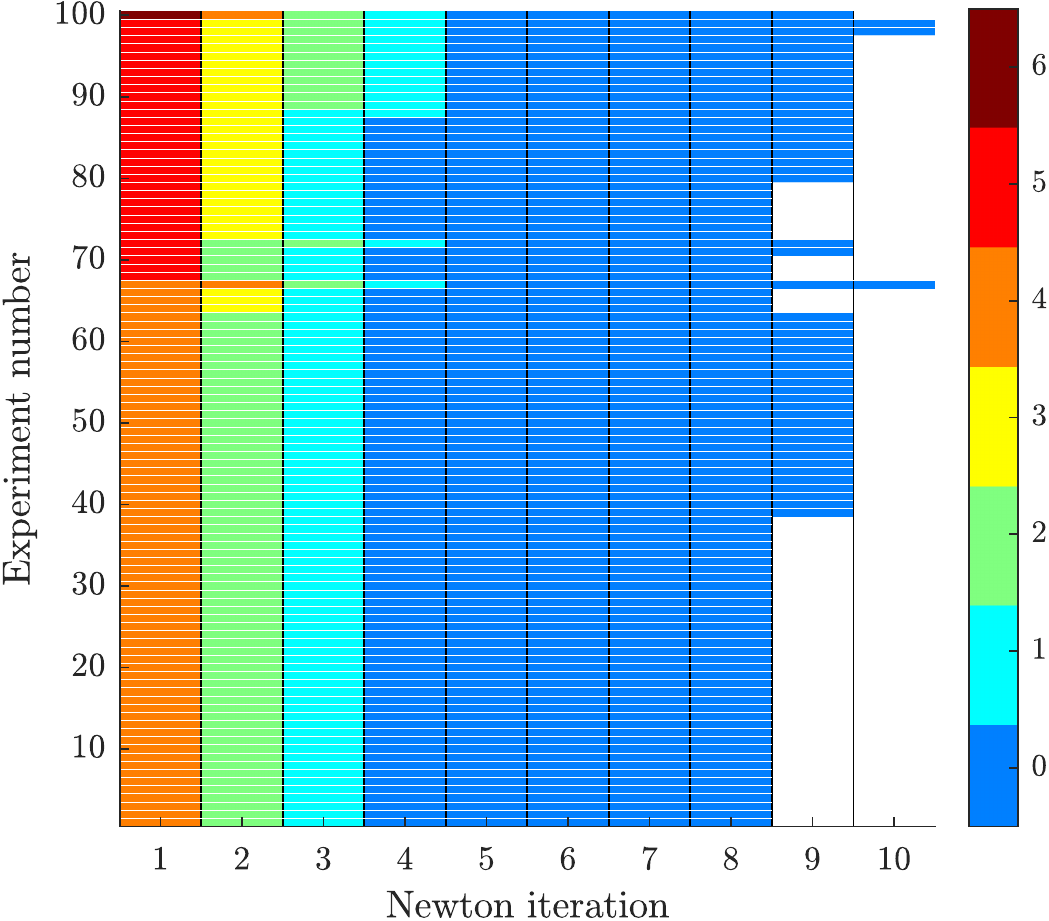}
\end{center}
\caption{\label{fig:backtrack}
The colour bar shows the number of backtracking steps for each Newton iteration, for each of 100 random log-normal microstructures with 100,000 grains (see Example \ref{example:back tracking}). Each row in the figure is a single draw, the number of coloured blocks indicates the number of Newton iterations until the error tolerance is reached, and the colour of each block indicates the number of backtracking steps per Newton iteration. The draws have been sorted by the number of backtracking steps in each Newton iteration.}
\end{figure}
\end{example}     
     
%\medskip     

%\clearpage
     
%\paragraph{Acknowledgements} 
\noindent
\textit{Acknowledgements.}
The authors thank Piet Kok, Wil Spanjer %, 
and
Carola Celada-Casero
%and Karo Sedighiani 
for fruitful discussions.
DPB acknowledges financial support from the EPSRC grant EP/V00204X/1. MP thanks the Centre for Doctoral Training in Mathematical Modelling, Analysis and Computation %(MAC-MIGS), 
funded by EPSRC grant 
%number 
EP/S023291/1.   
     
%\vspace{-0.5cm}     

\bibliographystyle{plain}
\bibliography{Bibliography}
   
\end{document}